\documentclass[english]{amsart}

\usepackage[francais,english]{babel}
\usepackage{amsthm}
\usepackage{latexsym}

\usepackage{amsmath}
\usepackage{amsfonts}
\usepackage{amssymb}
\usepackage{verbatim}
\usepackage{psfrag}
\usepackage{graphicx}
\newcommand{\drawat}[3]{\makebox[0pt][l]{\raisebox{#2}{\hspace*{#1}#3}}}

\newcommand{\Pn}{\mathbb{P}^2}
\newcommand{\Rn}{\mathbb{R}^2}

\theoremstyle{theorem}
\newtheorem{lemm}{Lemma}[subsection]
\newtheorem{prop}[lemm]{Proposition}
\newtheorem{coro}[lemm]{Corollary}
\newtheorem{conj}[lemm]{Conjecture}
\newtheorem{theo}{Theorem}

\newtheorem{defi}[lemm]{Definition}

\newtheorem{rema}[lemm]{Remark}

\title[Best polynomial bounds]{The best polynomial bounds for the number of triangles in a simple arrangement of $n$ pseudo-lines}
\author{J\'er\'emy Blanc}
\thanks{The author acknowledge support from the Swiss national science foundation}
\address{J\'er\'emy Blanc,
Universit\'e de Grenoble I,
UFR de Math\'ematiques,
UMR 5582 du CNRS, 
Institut Fourier, BP 74, 
38402 Saint-Martin d'H\`eres, France}

\begin{document}
\begin{abstract}It is well-known that affine (respectively projective) simple arrangement of $n$ pseudo-lines may have at most $n(n-2)/3$ (respectively $n(n-1)/3$) triangles. However, these bounds are reached for only some values of $n \pmod{6}$. We provide the best polynomial bound for the affine and the projective case, and for each value of $n \pmod{6}$.\end{abstract}
\maketitle
\section{Introduction}
Let us remind the reader some classical notions.
A \emph{pseudo-line} $L$ is a curve of the plane (which may be either the affine real plane $\mathbb{R}^2$ or the projective real plane $\Pn=\Pn(\mathbb{R})$) that may be stretched to a line, i.e. such that there exists an homeomorphism of the plane that sends $L$ on a line.

A \emph{simple affine (respectively projective) arrangement $\mathcal{A}$ of $n$ pseudo-lines} is a set $\mathcal{A}=\{A_1,..,A_n\}$ of $n$ pseudo-lines of the affine plane (respectively of the projective plane) such that 
two of the curves intersect transversally into exactly one point, and this one does not belong to any other curve of $\mathcal{A}$. In the sequel, $n$ will always be assume to be at least equal to $3$. There are exactly $n-1$ intersections on each pseudo-line, that delimit $n-2$ bounded segments and $2$ unbounded segments (respectively $n-1$ segments).


A particular case is when each pseudo-line is a line, i.e. when $\mathcal{A}$ is a \emph{simple arrangement of lines}. An arrangement of pseudo-lines is said to be \emph{stretchable} if there exists an homeomorphism of the plane that sends each pseudo-line on a line (such a transformation exists for each pseudo-line, here we ask that the same homeomorphism suits for each pseudo-line). There exist arrangements of pseudo-lines that are not stretchable (see \cite{bib:Gru}); the notion of arrangements of pseudo-lines is thus a non--trivial generalisation of those of arrangements of lines. Note that usually it is very hard to decide whether an arrangement is stretchable or not. 

A simple affine arrangement $\mathcal{A}$ of $n$ pseudo-lines decomposes the plane $\Rn$ into $n(n+1)/2+1$ regions, where $n(n-3)/2+1$ are bounded and $2n$ are unbounded. We say that these regions are the \emph{polygons} of the arrangement $\mathcal{A}$ and denote by $a_3(\mathcal{A})$ the number of bounded triangles (which are polygons delimited by three bounded segments). Since a bounded segment may not delimit two different triangles, and the number of bounded segments is $n(n-2)$, we have $a_3(\mathcal{A})\leq n(n-2)/3$.

The projective case is similar; a simple projective arrangement $\mathcal{A}$ of $n$ pseudo-lines decomposes $\Pn$ into ${n(n-1)}/2+1$ regions (once again called the \emph{polygons} of $\mathcal{A}$), and we denote by $p_3(\mathcal{A})$ the number of regions that are triangles. The number of segments is now $n(n-1)$, and if $n\geq 4$ a segment does not delimit two different triangles, thus we have $p_3(\mathcal{A})\leq n(n-1)/3$.

As we will prove later, these two rough bounds are reached for infinitely many values of $n$ (see Theorems~\ref{thm:AffBound}  and~\ref{thm:ProjBound} and also \cite{bib:Ha1}, \cite{bib:Ha2}, \cite{bib:Rou1}). But to be reached, the bound has to be an integer, and this implies conditions on the values of $n \pmod 3$. Furthermore, the affine bound is never reached if $n$ is even (Theorem~\ref{thm:AffBound}, \cite[Theorem 1.1]{bib:BaBlLo}) and the projective bound is never reached if $n$ is odd (Theorem~\ref{thm:ProjBound}, \cite[Theorem 2.21, p.26]{bib:Gru}). The number of triangles depends thus a lot from the value of $n\pmod 6$.

 
 \bigskip
 
 Denoting by $\overline{p_3^s}(n)$ (respectively $\overline{a_3^s}(n)$) the maximal number of triangles in a simple projective (respectively affine) arrangement of $n$ pseudo-lines, \emph{we look for the best polynomial upper bounds} for the values of $\overline{p_3^s}(n)$ and $\overline{a_3^s}(n)$, for each value of $n\pmod 6$. Furthermore, we would like to describe the values were the bounds are not reached (which seem to be rare). 

In the affine case, for $n\equiv 3,5\pmod{6}$, \cite[Theorem 1]{bib:Ha2} shows that the rough bound ($n(n-2)/3$) is the good one to choose, as this one is reached for infinitely many values. For $n\equiv 1\pmod{6}$, it was proved in \cite[Theorem 1.3]{bib:BaBlLo} that $\lfloor n(n-2)/3\rfloor=(n(n-2)-5)/3$ is the right bound. For $n$ even, the rough bound is not tight at all; it was proved in \cite[Theorem 1.1]{bib:BaBlLo} that $n(n-7/3)/3$ is an upper bound. However, this one is also not tight. We will improve this bound by showing in Corollary~\ref{Coro:nn5/2} that $n(n-5/2)/3$ suits, and that this one is tight (after rounding down for $n\equiv 2\pmod{6}$), since it is reached for infinitely many values of $n$.

In the projective case, for $n\equiv 0,4\pmod{6}$, the rough bound ($n(n-1)/3$) is tight, as it was proved in \cite{bib:Ha1}, \cite{bib:Ha2} and \cite{bib:Rou1}; more examples that reach the bound (in particular for any $n\leq 40$) were given in \cite{bib:BoRoSt}. For $n\equiv 2\pmod{6}$, it was observed in \cite[Theorem 3.2]{bib:Rou1} (see also \cite[Proposition 1.2]{bib:BaBlLo}) that $\lfloor n(n-1)/3\rfloor$ is never reached, i.e. that $(n(n-1)-5)/3$ is an upper bound. However, there was -- up to now -- no known example that reached this bound. We will provide an example of $26$ pseudo-lines in Figure~\ref{25pseudo} and show that the bound is reached for infinitely many values. Furthermore, we will show that the value of $p_3^s(n)$ for $n=8,14,20$, is $(n(n-1)-8)/3$, i.e. that the bound is not reached  for $n<26$. If $n$ is odd, it was observed by Granham (see \cite[page 26, Theorem 2.21]{bib:Gru}) that $p_3^s(n)\leq n(n-2)/3$. Reducing to the affine case, this bound is clearly tight (after rounding for $n\equiv 1\pmod{6}$).

Summing up, we will prove the following theorems:
\begin{theo}[polynomial bounds for affine arrangements]
\label{thm:AffBound}
Let $\mathcal{A}$ an affine arrangement of $n$ pseudo-lines. Then
$$a_3(\mathcal{A})\leq \left\{\begin{array}{llll}
n(n-5/2)/3&\mbox{ if }&n\equiv 0,4&\pmod{6}\\
(n(n-2)-2)/3&\mbox{ if }&n\equiv 1&\pmod{6}\\
(n(n-5/2)-2)/3&\mbox{ if }&n\equiv 2&\pmod{6}\\
n(n-2)/3&\mbox{ if }&n\equiv 3,5&\pmod{6}\end{array}\right.$$
Furthermore, for any integer $0\leq k\leq 5$, there exists infinitely many integers $n$ such that the bound above is reached, for $n\equiv k\pmod{6}$.
\end{theo} 

\begin{theo}[polynomial bounds for projective arrangements]
\label{thm:ProjBound}
Let $\mathcal{A}$ a projective arrangement of $n$ pseudo-lines, with $n\geq 4$. Then
$$p_3(\mathcal{A})\leq \left\{\begin{array}{llll}n(n-1)/3&\mbox{ if }&n\equiv 0,4&\pmod{6}\\
(n(n-2)-2)/3&\mbox{ if }&n\equiv 1&\pmod{6}\\
(n(n-1)-5)/3&\mbox{ if }&n\equiv 2&\pmod{6}\\
n(n-2)/3&\mbox{ if }&n\equiv 3,5&\pmod{6}\end{array}\right.$$
Furthermore, for any integer $0\leq k\leq 5$, there exists infinitely many integers $n$ such that the bound above is reached, for $n\equiv k\pmod{6}$.
\end{theo}
Furthermore, we will compute the explicit values of $\overline{p_3^s}(n)$ and $\overline{a_3^s}(n)$ for any integer $n\leq 30$:
\begin{theo}
\label{thm:30}
The bound of Theorem~\ref{thm:AffBound} is reached for any integer $n\leq 30, n\not=11,12$; and we have $\overline{a_3^s}(11)=32$ and $\overline{a_3^s}(12)=37$.

The bound of Theorem~\ref{thm:ProjBound} is reached for any integer $n\leq 30$, except for $n=8,11,12,14,20$, where the value of $\overline{p_3^s}(n)$ is respectively $16, 32, 40, 58, 124$.
\end{theo}
Furthermore, we conjecture the following statement (that was already stated in \cite[Conjecture 4.3]{bib:Rou2}, in the special case of projective arrangement of $n$ pseudo-lines, with $n\equiv 0,4\pmod{6}$)
\begin{conj}
The bounds of Theorems~\ref{thm:AffBound} and~\ref{thm:ProjBound} are reached for any integer $n\geq 21$.
\end{conj}
The smallest values of $n$ where the conjecture is still open are $31, 32$, $37$, $38$, $43$, $44$, $47$, $48$, $55$, $56$,...

\bigskip 

Our choice, in this article, is to be as comprehensive and self-contained as possible. To this aim, we will show provide some short proofs of some results cited above, in Lemmas~\ref{Lem:SegmPerdProjOdd},~\ref{Lemm:Bound2mod6Proj} and in Corollary~\ref{Coro:oddProj}, and will be able to prove the theorems  directly without any result.

This article is divided in the following way. Section \ref{Sec:ImportParity} deals with the parity of the number of pseudo-lines of the arrangements, and explains why this is so important. In particular, we prove (Corollary \ref{Coro:nn5/2}) that $\overline{a_3^s}(n)\leq n(n-5/2)/3$ if $n$ is even, a new tight bound stated before. We also provide the proof of the well-known inequality $\overline{p_3^s}(n)\leq n(n-2)/3$ for $n$ odd.
In Section \ref{Sec:proj2mod6}, we study the projective arrangements of $n$ pseudo-lines with $n\equiv 2\pmod{6}$. We give a simple proof of the inequality $\overline{p_3^s}(n)\leq (n(n-2)-5)/3$, already proved in \cite{bib:Rou1}, and provide the first known arrangement that reaches the bound (Figure \ref{25pseudo}).
In Section \ref{Sec:Doubl}, we give an adaptation of the doubling method of \cite{bib:Ha1} and \cite{bib:Rou1} to provide arrangements of $2n-1$ (respectively $2n-2$) pseudo-lines reaching the bounds of the theorems, starting with arrangements of $n$ pseudo-lines reaching the bound, where $n$ is odd (respectively even).
Finally, in Section \ref{Sec:Explicit}, we provide explicit arrangements that reach the bound (which have already appear in some previous works), and give the proofs of the Theorems.

The author would like to express his gratitude to J.P. Roudneff for his comments and remarks on the present article.



\section{The importance of the parity of the number of pseudo-lines}\label{Sec:ImportParity}
\begin{defi}
Let $\mathcal{A}$ be an affine or projective simple arrangement, and let $s$ be a segment of this arrangement. We say that $s$ is \emph{used} if it is an edge of a triangle of $\mathcal{A}$, otherwise we say that it is \emph{unused}.
\end{defi}
\begin{rema}
In an affine simple arrangement, any unbounded segment is unused.
\end{rema}

The following Proposition will explain why an affine arrangement with an even number of pseudo-lines contains many unused segments.

\begin{prop}\label{Prp:SegmPerd}
Let $\mathcal{A}$ be an \emph{affine} arrangement of $n$ pseudo-lines, where $n$ is even and $n\geq 4$.
Then, for any pseudo-line $L\in \mathcal{A}$, there exists an \emph{unused} bounded segment $s\subset L'\in\mathcal{A}$ whose intersection with  $L$ consists of one point (in particular $s \not\subset L$).
\end{prop}
\begin{proof}
Let $L\in\mathcal{A}$ be a pseudo-line. Applying an homeomorphism of the plane to $\mathcal{A}$, we may assume that $L$ is the line $y=0$ and that the points of intersection of $L$ with the other lines of $\mathcal{A}$ are $(1,0), (2,0),..., (n-1,0)$.
For $i=1,...,n-1$, we denote by $R_i$ the pseudo-line of $\mathcal{A}\backslash\{L\}$ that passes through the point $(i,0)$ and by $r_i^{+}$ (respectively $r_i^{-}$) the segment of $R_i$ that touches $(i,0)$ and belongs to the upper (respectively lower) half-plane. For $i=1,...,n-2$,
we call $l_i$ the segment  $\{ (x,0)\in\Rn \ |\ i\leq x \leq i+1\}$ of the line $L$; we call also $l_0=\{ (x,0) \in \Rn\ |\ x<1\}$ and $l_{n-1}=\{ (x,0) \in \Rn\ |\ x>n-1\}$, the two unbounded segments of $L$. 
For $i=0,...,n-1$, we denote by $p_i^{+}$ (respectively $p_i^{-}$) the polygon of $\mathcal{A}$ -- that may be bounded or not -- that touches $l_i$, and belongs to the upper (respectively lower) half-plane. 
 
\begin{figure}[ht]{\begin{center}\includegraphics[width=10cm]{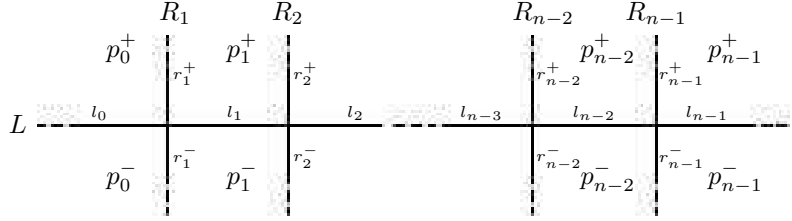}
\drawat{-85mm}{26mm}{$R_1$}%
\drawat{-105mm}{11mm}{$L$}%
\drawat{-70mm}{26mm}{$R_2$}%
\drawat{-38mm}{26mm}{$R_{n-2}$}%
\drawat{-23mm}{26mm}{$R_{n-1}$}%
\drawat{-94mm}{13.2mm}{\tiny $l_0$}%
\drawat{-76mm}{13.2mm}{\tiny $l_1$}%
\drawat{-76mm}{21.2mm}{$p_1^{+}$}%
\drawat{-76mm}{4.2mm}{$p_1^{-}$}%
\drawat{-12mm}{21.2mm}{$p_{n-1}^{+}$}%
\drawat{-12mm}{4.2mm}{$p_{n-1}^{-}$}%
\drawat{-29mm}{21.2mm}{$p_{n-2}^{+}$}%
\drawat{-29mm}{4.2mm}{$p_{n-2}^{-}$}%
\drawat{-92mm}{21.2mm}{$p_0^{+}$}%
\drawat{-92mm}{4.2mm}{$p_0^{-}$}%
\drawat{-60mm}{13.2mm}{\tiny $l_2$}%
\drawat{-45mm}{13.2mm}{\tiny $l_{n-3}$}%
\drawat{-30mm}{13.2mm}{\tiny $l_{n-2}$}%
\drawat{-15mm}{13.2mm}{\tiny $l_{n-1}$}%
\drawat{-83.2mm}{18.2mm}{\tiny $r_{1}^{+}$}%
\drawat{-83.2mm}{7.2mm}{\tiny $r_{1}^{-}$}%
\drawat{-67.2mm}{18.2mm}{\tiny $r_{2}^{+}$}%
\drawat{-67.2mm}{7.2mm}{\tiny $r_{2}^{-}$}%
\drawat{-34.9mm}{18.2mm}{\tiny $r_{n-2}^{+}$}%
\drawat{-34.9mm}{7.2mm}{\tiny $r_{n-2}^{-}$}%
\drawat{-18.6mm}{18.2mm}{\tiny $r_{n-1}^{+}$}%
\drawat{-18.6mm}{7.2mm}{\tiny $r_{n-1}^{-}$}%
\end{center}}
\caption{The arrangement $\mathcal{A}$, in a neighbourhood of the pseudo-line $L$.\label{FigurePropSegmPerd}}\end{figure}

 \begin{enumerate}
 \item
{\it Assume that $p_1^{+}$ and $p_{n-2}^{-}$ are triangles.}\\
This implies that neither $r_1^{-}$ nor $r_{n-1}^{+}$ are used. If  $L_1$ and $L_n$ intersect in the upper (respectively lower) half-plane, then $r_{n-1}^{+}$ (respectively $r_{1}^{-}$) is bounded. 
\item
{\it Assume that neither $p_1^{+}$ nor $p_1^{-}$ is a triangle.}\\
In this case, neither $r_1^{+}$ nor $r_1^{-}$ is a used segment of $\mathcal{A}$. Since one of these is bounded, we are done.
\item
We may now assume that $p_1^{+}$ is a triangle (the case where $p_1^{-1}$ is a triangle is the same). Suppose, for contradiction, that every segment $r_i^{\pm}$ that is bounded is used, for $i=1,...,n-1$. \\
{\it Claim:} for $i\in\{1,..,n-1\}$ even (respectively odd) $p_i^{+}$ (respectively $p_i^{-}$) is not a triangle; and if neither $p_{i-1}^+$ nor $p_{i-1}^{-}$ is a triangle, then $p_{i}^{-}$ (respectively $p_{i}^{+}$) is a triangle.\\
The claim being true for $i=1$, we prove it by induction on $i$. Assume that $i$ is even (the odd case is similar, exchanging $+$ and $-$). If $p_{i-1}^{+}$ is a triangle, then $p_{i}^{+}$ is not a triangle, since two triangles do not have a common edge, and the claim is proved.  Assume that $p_{i-1}^{+}$ is not a triangle, which implies that $i>2$, and let us illustrate the situation (by induction hypothesis, neither $p_{i-1}^{-}$ nor $p_{i-2}^{+}$ is a triangle).
\begin{figure}[ht]{\begin{center}\includegraphics[width=11cm]{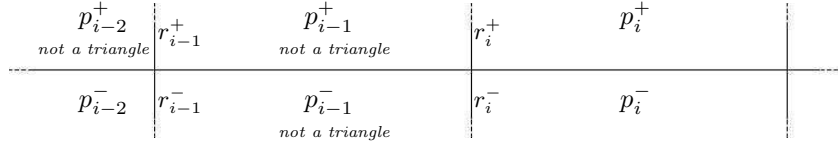}
\drawat{-49.5mm}{13mm}{\small $r^{+}_{i}$}%
\drawat{-49.5mm}{4mm}{\small $r^{-}_{i}$}%
\drawat{-91.5mm}{13mm}{\small $r^{+}_{i-1}$}%
\drawat{-91.5mm}{4mm}{\small $r^{-}_{i-1}$}%
\drawat{-30mm}{15mm}{$p_{i}^{+}$}%
\drawat{-30mm}{4mm}{$p_{i}^{-}$}%
\drawat{-72mm}{15mm}{$p_{i-1}^{+}$}%
\drawat{-75.5mm}{0mm}{\tiny{\it not a triangle}}%
\drawat{-75.5mm}{11mm}{\tiny{\it not a triangle}}%
\drawat{-72mm}{4mm}{$p_{i-1}^{-}$}%
\drawat{-102mm}{15mm}{$p_{i-2}^{+}$}%
\drawat{-107.5mm}{11mm}{\tiny{\it not a triangle}}%
\drawat{-102mm}{4mm}{$p_{i-2}^{-}$}%
\end{center}} \caption{The situation of the polygons $p_{i-2}^{\pm},p_{i-1}^{\pm},p_i^{\pm}$\label{FigDrLemm2}.}\end{figure}

 The segment $r_{i-1}^{+}$ -- that is common to $p_{i-1}^{+}$ and $p_{i-2}^{+}$ -- is unused and consequently unbounded. This implies that $R_i$ and $R_{i-1}$ intersect into the lower half-plane, so $r_i^{-}$ is bounded. The assumption made above implies that $r_i^{-}$ is used, so $p_i^{-}$ is a triangle, and not $p_i^{+}$. This achieves to prove the claim.
 
 Applying the claim for $i=n-1$, the region $p_{n-1}^{-}$ is not a triangle. If $p_{n-1}^{+}$ is a triangle, we obtain the case $1$ treated above. Otherwise, we obtain the case $2$ treated above, after reversing the order of the segments on $L$.
\end{enumerate}
\end{proof}
This Proposition yields the following bound,  which ameliorates the one of \cite{bib:BaBlLo}, that was $n(n-7/3)/3$.
\begin{coro}\label{Coro:nn5/2}
Let $n$ be an even integer and let $\mathcal{A}$ be a simple affine arrangement of $n$ pseudo-lines. Then $a_3(\mathcal{A})\leq \lfloor n(n-5/2)/3\rfloor$.
\end{coro}
\begin{proof}
According to Proposition~\ref{Prp:SegmPerd}, we may associate to any pseudo-line $A\in \mathcal{A}$ a segment $s_A$ of $\mathcal{A}$, that is bounded, unused and whose intersection with $A$ consists of one point. Observe that a segment may be associated to at most two lines, which implies that the number of unused segments is at least $n/2$. The number of used segments is consequently at most $n(n-5/2)$, whence $a_3(\mathcal{A})\leq \lfloor n(n-5/2)/3\rfloor$.
\end{proof}

Similarly to the affine arrangements with an even number of pseudo-lines, any projective arrangement with an odd number of pseudo-lines contains many unused segments. This was firstly observed by Granham (see \cite[page 26, Theorem 2.21]{bib:Gru}) and used many times, we prove once again this simple result:
\begin{lemm}\label{Lem:SegmPerdProjOdd}
Let $\mathcal{A}$ be a projective arrangement of $n$ pseudo-lines, where $n\geq 5$ is an odd integer. Then, any pseudo-line of $\mathcal{A}$ contains at most $n-2$ used segments, i.e. there exists at least one unused segment on each pseudo-line.
\end{lemm}
\begin{proof}
Let $L\subset \mathcal{A}$ be a pseudo-line of $\mathcal{A}$.
Then, $L$ contains exactly $n-1$ segments and touches exactly $2(n-1)$ polygons of the arrangement.
The lemma follows from the following two observations:
\begin{enumerate}
\item
Since $n$ is odd, any set of $n-1$ polygons touching $L$ contains at least two polygons that have a common edge (this is false if $n$ is even).
\item
Since the number of pseudo-lines is at least $5$, two distinct triangles of $\mathcal{A}$ do not have a common edge. 
\end{enumerate}
\end{proof}
This Lemma implies -- as it was already noticed in \cite{bib:Gru} -- the following bound.
\begin{coro}\label{Coro:oddProj}
Let $n$ be an odd integer and let $\mathcal{A}$ be a simple projective arrangement of $n$ pseudo-lines. Then $p_3(\mathcal{A})\leq \lfloor n(n-2)/3\rfloor$.
\end{coro}
\begin{proof}
According to Lemma~\ref{Lem:SegmPerdProjOdd}, there is an unused segment on each pseudo-line of $\mathcal{A}$. The number of used segments is thus at most $n(n-2)$, whence $p_3(\mathcal{A})\leq \lfloor n(n-2)/3\rfloor$.\end{proof}


\section{Projective arrangements of $n$ pseudo-lines, with $n\equiv 2 \pmod{6}$}\label{Sec:proj2mod6}
The most interesting case in our subject is when we have a projective simple arrangement $\mathcal{A}$ of $n$ pseudo-lines, where $n\equiv 2\pmod{6}$. The rough bound on the number of segment implies that $$p_3(\mathcal{A})\leq \lfloor n(n-2)/3\rfloor=(n(n-2)-2)/3.$$
However, as it was already observed in \cite[Theorem 3.2]{bib:Rou1} and \cite[Proposition 1.2]{bib:BaBlLo}, the above bound is never reached. We prove once again this result for self-containedness.
\begin{lemm}\label{Lemm:Bound2mod6Proj}
Let $\mathcal{A}$ be a projective arrangement of $n$ pseudo-lines with $n\equiv 2\pmod{6}$, then $p_3(\mathcal{A})\leq (n(n-2)-5)/3$. Furthermore, if the equality occurs, the $5$ unused segments of $\mathcal{A}$ delimit a pentagon of the arrangement.
\end{lemm}
\begin{proof}
Since there are $n(n-1)$ segments, and $n(n-1)\equiv 2\pmod{3}$, there exists at least two segments of $\mathcal{A}$ that are unused.
Let $A\in \mathcal{A}$ be a pseudo-line and let $s$ be an unused segment contained in $A$, whose extremities are the points $v_1$ and $v_2$. There are $6$ other segments, and $6$ polygons that touch either $v_1$ or $v_2$. We call these respectively $s_1,...,s_6$ and $p_1,...,p_6$ as in Figure~\ref{Fig52}.

\begin{figure}[ht]{\begin{center}\includegraphics[width=5cm]{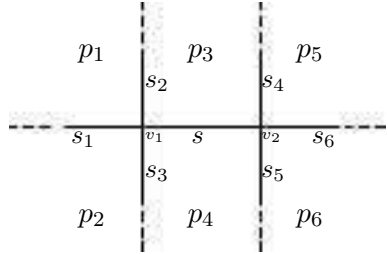}
\drawat{-27.5mm}{26mm}{\large $p_{3}$}%
\drawat{-27.5mm}{4mm}{\large $p_{4}$}%
\drawat{-13mm}{26mm}{\large $p_{5}$}%
\drawat{-42mm}{4mm}{\large $p_{2}$}%
\drawat{-42mm}{26mm}{\large $p_{1}$}%
\drawat{-13mm}{4mm}{\large $p_{6}$}%
\drawat{-11mm}{14.5mm}{$s_6$}%
\drawat{-27mm}{14.5mm}{$s$}%
\drawat{-43mm}{14.5mm}{$s_1$}%
\drawat{-33.2mm}{15mm}{\tiny $v_1$}%
\drawat{-17.7mm}{15mm}{\tiny $v_2$}%
\drawat{-33.2mm}{22mm}{$s_2$}%
\drawat{-33.2mm}{10mm}{$s_3$}%
\drawat{-17.7mm}{22mm}{$s_4$}%
\drawat{-17.7mm}{10mm}{$s_5$}%
\end{center}} \caption{The situation of the polygons $p_1,...,p_6$ and the segments $s,s_1,...,s_6$\label{Fig52}.}\end{figure}

Since $s$ is unused, neither $p_3$ nor $p_4$ is a triangle. The number of pseudo-lines being at least $3$, two triangles do not have a common edge. Thus, either $p_1$ or $p_2$ is not a triangle, which implies that either $s_2$ or $s_3$ is unused. Similarly, either $s_4$ or $s_5$ is unused. 

This shows that each unused segment of $\mathcal{A}$ touches two other unused segments. The number of unused segments being equal to $2\pmod{3}$, there are at least $5$ unused segments and if there are exactly $5$ of them, these form a pentagon of the arrangement.
\end{proof}

\begin{figure}[ht]{\begin{center}\includegraphics[width=11.5cm]{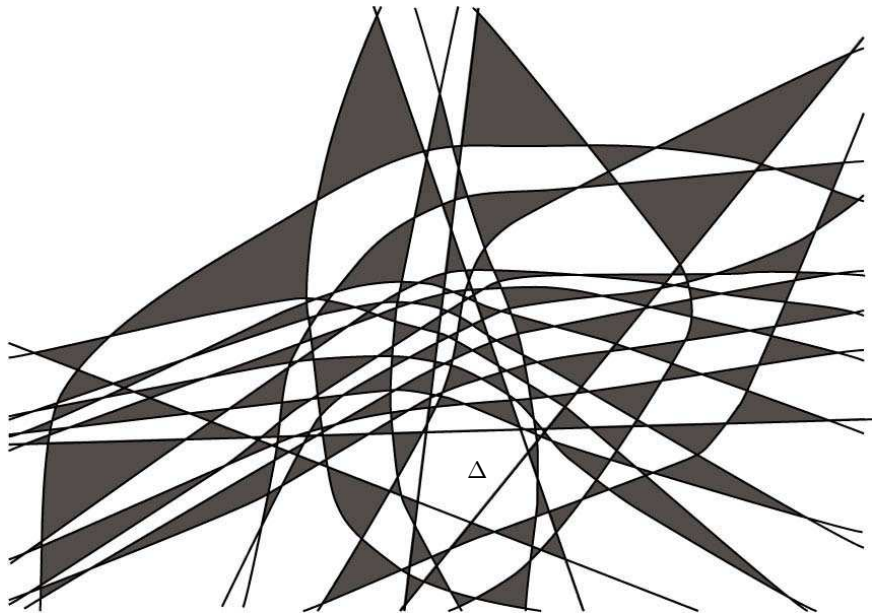}
\drawat{-55mm}{18mm}{\small $\Delta $}%
\end{center}}
\caption{An affine arrangement that gives, with the line at infinity a maximal projective arrangement of $n=26$ pseudo-lines with $(n(n-1)-5)/3=215$ triangles, the maximum possible. There are exactly five unused segments, that form a pentagon $\Delta$.\label{25pseudo}}\end{figure}

Up to now, there was no known example of arrangement that reached this bound. We provide in Figure~\ref{25pseudo} the first example that reaches the bound, which is an arrangement of $26$ pseudo-lines. Furthermore, an exhaustive counting (that may be at hand for $14$ pseudo-lines and with computer for $20$ pseudo-lines) shows that no such example exists for $8$, $14$ and $20$ pseudo-lines.
\section{The doubling method}\label{Sec:Doubl}
Given a projective arrangement of $n$ pseudo-lines with $n(n-1)/3$ triangles -- which may occur only if $n\equiv 0,4 \pmod{6}$ -- \cite{bib:Ha1} and \cite{bib:Rou1} give a way to construct a projective arrangement of $2n-2$ pseudo-lines with $(2n-2)(2n-3)/3$ triangles. We generalise this construction to arrangements that may be affine or projective and having any number of pseudo-lines (a similar construction for affine arrangements of lines may be found in \cite{bib:BaBlLo}).
\begin{lemm}
Let $\mathcal{A}$ be an affine (respectively projective) arrangement of $n$ pseudo-lines, where $n$ is odd (respectively even), such that one pseudo-line of $\mathcal{A}$ touches $n-2$ (respectively $n-1$) triangles of the arrangement.
Then, there exists an affine (respectively projective) arrangement of $2n-1$ (respectively $2n-2$) pseudo-lines $\mathcal{B}$, such that $a_3(\mathcal{B})=a_3(\mathcal{A})+(n-1)^2$ (respectively $p_3(\mathcal{B})=p_3(\mathcal{A})+(n-1)(n-2)$).
\end{lemm}
\begin{proof}
If the arrangement is projective we choose one pseudo-line, we stretch it and put it at infinity, to obtain an affine arrangement having one pseudo-line less than the projective arrangement. 
 
We work with an affine arrangement $\mathcal{A}=\{A_1,...,A_{m}\}$ of $m$ pseudo-lines, where $m$ is an odd integer and $A_m$ touches exactly $m-2$ triangles of the arrangement, and denote by $\overline{\mathcal{A}}=\mathcal{A}\cup\{L_{\infty}\}$ the projective arrangement obtained by adding to $\mathcal{A}$ the line at infinity. We will construct $m$ pseudo-lines $B_1,...,B_m$ such that the affine arrangement $\mathcal{B}=\{A_1,...,A_{m-1},B_1,...,B_m\}$ has $(m-1)^2$ more triangles than $\mathcal{A}$, and such that  if $A_m$ touches exactly $m$ triangles of the projective arrangement $\overline{\mathcal{A}}$ then the projective arrangement $\overline{\mathcal{B}}=\mathcal{B}\cup\{L_{\infty}\}$ has $m(m-1)$ triangles more than $\overline{\mathcal{A}}$. This will achieve the proof of the Lemma.

By an homeomorphism of the affine plane, we may assume that $A_m$ is the horizontal line $y=0$, and up to a reordering, that the $x$-coordinate of $A_m \cap A_i$ is lower than the $x$-coordinate of $A_m \cap A_{i+1}$ for $i=1,...,m-2$. We write $A_0=L_{\infty}$ and for $i=0,...,m-1$ we denote by $s_i$ the union of the two segments of the arrangement $\overline{\mathcal{A}}$ that lie on $A_i$ and touch $A_m$; this yields the situation described in Figure~\ref{dupli1} (after assuming that $s_1$ and $s_2$ meet in the upper half-plane, which implies that $s_{m-2}$ and $s_{m-1}$ do the same, since $m$ is odd). In particular, $s_i$ intersects $s_{i+1}$ for $i=1,...,m-2$.

\begin{figure}[ht]{\begin{center}\includegraphics[width=10cm]{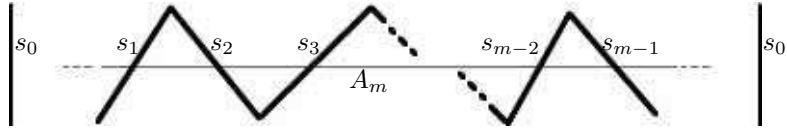}
\drawat{-87mm}{10mm}{$s_1$}%
\drawat{-74.5mm}{10mm}{$s_2$}%
\drawat{-63mm}{10mm}{$s_3$}%
\drawat{-56mm}{5mm}{$A_m$}%
\drawat{-38.5mm}{10mm}{$s_{m-2}$}%
\drawat{-22.5mm}{10mm}{$s_{m-1}$}%
\drawat{-1mm}{10mm}{$s_{0}$}%
\drawat{-100.5mm}{10mm}{$s_{0}$}%
\end{center}} \caption{The sequence of the segments $s_0,s_1,...,s_{m-1}$, near the line $A_0$\label{dupli1}.}\end{figure}

We construct the pseudo-lines of the set  $\mathcal{B}^{+}=\{B_1,..,B_m\}$ near $A_0$ in the following way:
\begin{enumerate}
\item
For $j=0,...,m-1$, each pseudo-line of $\mathcal{B}^{+}$ intersect $A_j$ on the segment $s_j$.
\item
For $j=0,...,m-1$, between the two segments $s_j$ and $s_{k}$ -- where $k=j+1$ if $j<m-1$ and $k=0$ otherwise -- one special pseudo-line in $\mathcal{B}^{+}$ does not touch any other curve of $\mathcal{B}^{+}$; otherwise each element of $\mathcal{B}^{+}$ intersect exactly one other. The special pseudo-line is the one whose intersection with $s_j$ has the lowest (respectively biggest) $x$-coordinate if $j$ is even (respectively if $j$ is odd).
\end{enumerate}
This yields the situation described in Figure~\ref{dupli2}.

\begin{figure}[ht]{\begin{center}\includegraphics[width=10cm]{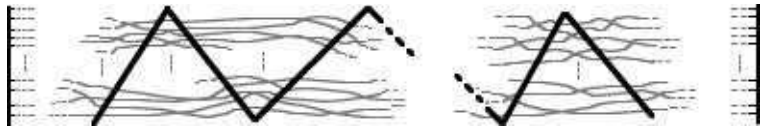}
\end{center}} \caption{The arrangement $\overline{\mathcal{B}}=\{A_0,A_1,...,A_{m-1},B_1,...,B_m\}$\label{dupli2}.}\end{figure}

Let us prove that this construction is possible, i.e. that the vertical ordering of the elements of $\mathcal{B}^{+}$ at the left of the affine part is the opposite of the one at the right, so that the ordering on the line at infinity is the same on each way. We order the curves so that $B_i$ is upper that $B_{i+1}$ at the left of the affine part, for $i=1,..,m-1$. Between two consecutive segments $s_j$ and $s_{k}$, we say that a pseudo-line of $\mathcal{B}^{+}$ \emph{is going down} (respectively \emph{is going up}) if it crosses a pseudo-line which is below it on the segment $s_j$ (we do as if the segments $s_j$ was almost vertical); we say that the pseudo-line that crosses no other element of $\mathcal{B}^{+}$ is \emph{horizontal}.
Observe that, from left to right, each pseudo-line $B_1,B_3,...,B_ {m-2}$ (respectively $B_2,B_4,...,B_{m-1}$) goes down (respectively goes up) until it is at the top (respectively at the bottom), stay horizontal, and then goes down (respectively goes up). The pseudo-line $B_m$ stay horizontal, and then goes up. Since there are $m$ places to cross, the line $B_i$ is at the place $m+1-i$ at the end.

Furthermore, each pseudo-line of $\mathcal{B}^{+}$ stays horizontal exactly at one place. This shows that each one intersects $m-1$ other pseudo-lines of $\mathcal{B}^{+}$.  The ordering at the left and the right being reversed, each two elements of $\mathcal{B}^{+}$ intersect into exactly one point.

The two arrangements $\mathcal{B}$ and $\overline{\mathcal{B}}$ are thus arrangement of pseudo-lines. Counting the new triangles added (for example by observing that no new unused segment is added), we obtain the results needed.
\end{proof}

\begin{coro}\label{coro:SeqBound}
If the bound of Theorem~\ref{thm:AffBound} is reached for some odd (respectively even) integer $n$, then it is reached for $2n-1$ (respectively for $2n-2$). The same is true for the bound of Theorem~\ref{thm:ProjBound}.
\end{coro}
\begin{proof} This may be proved by direct counting, or by observing that the number of unused segments of the "small" arrangement is the same as the number of unused segments of the arrangement obtained by the doubling method.
\end{proof}

\section{Explicit configurations that reach the bounds}\label{Sec:Explicit}
\begin{figure}[ht]{
\begin{center}
\includegraphics[width=12cm]{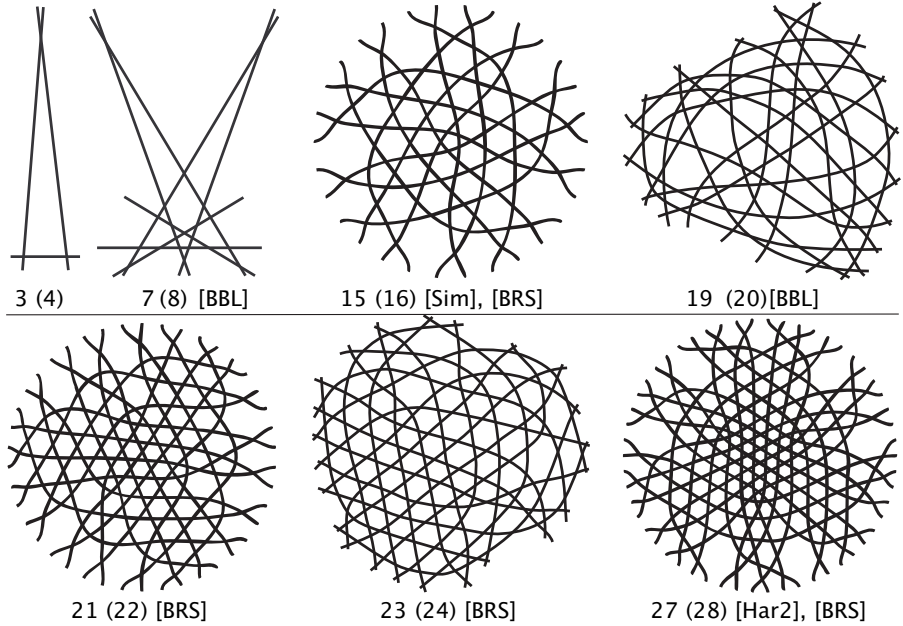} 
\end{center}}
\caption{Affine arrangements of $3,7,15,19,21,23,27$ pseudo-lines that reach the bounds of Theorems and that give -- adding a line far away --  affine (respectively projective) arrangements of $4,8,16,20,22,24,28$ (respectively $4,16,22,24,28$) pseudo-lines that reach the bounds.\label{FigureAffineBeg}}\end{figure}

The bounds of Theorems~\ref{thm:AffBound} and~\ref{thm:ProjBound} have been proved in Corollaries~\ref{Coro:nn5/2} and~\ref{Coro:oddProj}, and in Lemma~\ref{Lemm:Bound2mod6Proj}. To prove the two theorems, it remains to provide infinite sequences that reach the bound, for any value of $n\pmod{6}$ and for both affine and projective cases. This will be done by providing concrete examples (Figure~\ref{FigureAffineBeg}) and using Corollary~\ref{coro:SeqBound}.

\begin{prop}\label{Prp:Starts}
For $n=3,4,7,8,15,16,19,20,21,22,23,24,26,27,28$, the bounds of Theorems~\ref{thm:AffBound} (affine arrangements) is reached. For the same values -- except $n=3,8,20$ -- the bound of Theorem~\ref{thm:ProjBound} (projective arrangements) is reached.
\end{prop}
\begin{rema}
Except for $n=26$, the Proposition follows from \cite{bib:Sim}, \cite{bib:Ha2}, \cite{bib:BoRoSt} and \cite{bib:BaBlLo}; the configurations of these articles are given in Figure~\ref{FigureAffineBeg}.
\end{rema}
\begin{proof}
The odd values are provided directly by the arrangements of Figure~\ref{FigureAffineBeg}. Adding a line far away, we obtain the even values, except $n=26$. The projective arrangement of $26$ pseudo-lines provided in Figure~\ref{25pseudo} reaches the bound of Theorem~\ref{thm:ProjBound}; choosing the line at infinity near one of the five lines of the arrangements that touch the pentagon $\Delta$, we obtain an affine arrangement of $26$ pseudo-lines that reaches the bound of Theorem~\ref{thm:AffBound}.
\end{proof}

\begin{coro}\label{Coro:ExplicitInfSeq}
For $n=m\cdot 2^t+1$, where $t\geq 0$ and $m=4,6,14,18,20,22,26$,  the bounds of Theorems~\ref{thm:AffBound} and~\ref{thm:ProjBound} are reached.

For $n=m\cdot 2^t+2$, where $t\geq 0$ and $m=4,14,20,22,24,26$,  the bounds of Theorems~\ref{thm:AffBound} and~\ref{thm:ProjBound} are reached.
\end{coro}
\begin{proof}
According to Corollary~\ref{coro:SeqBound}, if the bounds are reached for $n=m\cdot 2^t+1$ (respectively for $n=m\cdot 2^t+2$), these are reached for $n=m\cdot 2^{t+1}+1$ (respectively for $n=m\cdot 2^{t+1}+2$). This corollary follows then directly from Proposition~\ref{Prp:Starts}.
\end{proof}

\subsection{Proofs of theorems~\ref{thm:AffBound},~\ref{thm:ProjBound} and~\ref{thm:30}} We are now ready to give the proofs of the Theorems.
\begin{proof}[Proof of Theorems~\ref{thm:AffBound} and~\ref{thm:ProjBound}]The bounds of Theorems~\ref{thm:AffBound} and~\ref{thm:ProjBound} have been proved in Corollaries~\ref{Coro:nn5/2} and~\ref{Coro:oddProj}, and in Lemma~\ref{Lemm:Bound2mod6Proj}.
The infinite sequences of values that reach the bound are provided by Corollary~\ref{Coro:ExplicitInfSeq}.
\end{proof}
\begin{proof}[Proof of Theorem~\ref{thm:30}]
The values of $\overline{a_3^s}(n)$ and $\overline{p_3^s}(n)$ for $n=3,4$ are easy to compute. 
Applying Corollary \ref{Coro:ExplicitInfSeq} for $m=4$ we obtain the values for $n=5,6,9,10,17$ and $18$. Similarly, we see that $\overline{a_3^s}(n)$ and $\overline{p_3^s}(n)$ reach the bound for the following values of $n$:
\[\begin{array}{|l|l|}
\hline
m & n\\
\hline
4 & 5,6,9,10,17,18,33,34,...\\
\hline
6 & 7,13,25,49,97,193,...\\
\hline
14 & 15,16,29,30, 57,58,...\\
\hline
18 & 19, 37, 73, 145, 289,...\\
\hline
\end{array} \hspace{1 cm} \begin{array}{|l|l|}
\hline
m & n\\
\hline
20 & 21,22,41,42,81,82,...\\
\hline
22 & 23, 24, 45, 46, 89, 90, ... \\
\hline
24 & 26, 50, 98, 194, 386, ... \\
\hline
26 & 27, 28, 53, 54, 105, 106,...\\
\hline
\end{array}\]
In particular, for $n\leq 30$, both bounds are reached, except possibly for $n=8,11,12,14,20$. 

In the projective case, the bound (of Theorem \ref{thm:ProjBound}) is not reached, for any of these $5$ values. For $n=12$, a nice proof may be found in \cite[Theorem 4.2]{bib:Rou2}. For $n=8$, it is quite simple to see it directly by hand, and for $n=14$ it remains possible, with a little bit more work. In fact, a simple computer algorithm proves that the value of $\overline{p_3^s}(n)$ for $n=8,11,12,14,20$ is respectively $16, 32, 40, 58, 124$.

In the affine case, the bound (of Theorem~\ref{thm:AffBound}) is reached for $n=8$ and $n=20$ (Proposition \ref{Prp:Starts}), and then also for $n=14$, according to Corollary~\ref{coro:SeqBound}. The remaining two cases ($n=11,12$) do not reach the bound. In particular,  $\overline{a_3^s}(11)=32$ and $\overline{a_3^s}(12)=37$. This may also be proved by computer (see \cite{bib:BaBlLo}).
\end{proof}
\begin{rema}
In \cite[Table 1]{bib:Rou3}, it is stated (without proof) that $\overline{p_3^s}(12)=42$, a result that has been cited in \cite{bib:BaBlLo}. This should be a misprint (as it was confirmed to us by J.P. Roudneff in a private communication). A projective simple arrangement of $12$ pseudo-lines with $42$ triangles would have only $6$ unused segments, which would form an hexagon; by hand it is quite easy to observe that this is impossible. 
\end{rema}
\section{Remarks on arrangements of lines}
Since an an arrangement of lines is also an arrangement of pseudo-lines, all the bounds of Theorems \ref{thm:AffBound} and \ref{thm:ProjBound} are still valid for affine and projective arrangements of lines. Moreover, it is possible to stretch the first three examples of Figure \ref{FigureAffineBeg} (the first two are arrangements of lines and the third is stretched in \cite{bib:Sim}). The duplicate construction may also be applied to straight lines (see \cite{bib:FoR} and \cite{bib:BaBlLo}). To obtain the proof of Theorems \ref{thm:AffBound} and \ref{thm:ProjBound} for lines, it remains in fact to provide an arrangement of $n$ lines with $n\equiv 2 \pmod{6}$ having exactly $(n(n-1)-5)/3$ triangles (which could be done by stretching the arrangement of Figure \ref{25pseudo}). 

\end{document}